# ANALYSIS OF BINARY SPATIAL DATA BY QUASI-LIKELIHOOD ESTIMATING EQUATIONS

By Pei-Sheng Lin and Murray K. Clayton

*National Chung Cheng University and University of Wisconsin–Madison*

The goal of this paper is to describe the application of quasi-likelihood estimating equations for spatially correlated binary data. In this paper, a logistic function is used to model the marginal probability of binary responses in terms of parameters of interest. With mild assumptions on the correlations, the Leonov–Shiryaev formula combined with a comparison of characteristic functions can be used to establish asymptotic normality for linear combinations of the binary responses. The consistency and asymptotic normality for quasi-likelihood estimates can then be derived. By modeling spatial correlation with a variogram, we apply these asymptotic results to test independence of two spatially correlated binary outcomes and illustrate the concepts with a well-known example based on data from Lansing Woods. The comparison of generalized estimating equations and the proposed approach is also discussed.

**1. Introduction.** This paper was originally motivated by a question arising in forest ecology. Specifically, for each of a number of trees in a plantation, it is possible to determine the status of the tree (alive or dead) and whether a particular insect pest is present on the tree. The ecological question of interest is whether tree status is independent of the presence/absence of the pest. Due to biotic interaction such as mutualism or competition, observations recorded at locations usually exhibit spatial autocorrelation and this renders the traditional methods unsuitable. The question then becomes: how do we test for independence of two binary variables in the presence of spatial autocorrelation?

To abstract the problem, suppose that we have observations $\mathbf{Y} = (Y(\mathbf{s}_1), \ldots, Y(\mathbf{s}_N))^T$ that are binary responses drawn from a random field on an $m \times n$ grid with $N = mn$, and where $\mathbf{s}_i \in R^2$ denotes the location of the $i$th binary variable.

The mean response $\boldsymbol{\theta} = E(\mathbf{Y})$ is assumed to be associated with the measurements of explanatory variables $\mathbf{T}$ through a link function $h(\boldsymbol{\theta}) = \mathbf{T}\boldsymbol{\beta}$. The problem of detecting association then turns to the estimation of the parameters $\boldsymbol{\beta}$.

There is relatively little literature on the estimation of fixed effects $\boldsymbol{\beta}$ when observations are binary and spatially dependent. Albert and McShane [2] and Gotway and Stroup [11] mentioned the use of generalized estimating equations (GEE) in this setting. However, the GEE approach has been most fully developed for longitudinal data [15]. These methods do not directly apply to spatial data, and therefore there is a need to develop an estimating equation approach corresponding to the covariance structures typically found in spatial statistics, which are neither diagonal nor block diagonal.

To proceed, we will rely on quasi-likelihood (QL) ideas. The concept of quasi-likelihood functions was first introduced by Wedderburn [25] for observations from the exponential family of distributions with only mean and variance being specified. McCullagh (year?) and McCullagh and Nelder (year?) broadened the application to multivariate cases and discuss asymptotic properties. A potential obstacle in applying QL methods is that most of the developed methods were mainly based on the assumption of independent observations. Although QL functions can be defined for dependent data, McCullagh and Nelder [20] raised some concerns about the application to dependent data.

There are some examples in the literature of the application of QL estimating equations to correlated data, although most of these are focused on count data. Zeger [26] developed QL estimating equations for time series count data for a specific type of covariance function, and this approach was also applied by McShane, Albert and Palmatier [21] to spatial data. Papers more closely related to the current one are those of Heagerty and Lumley [13] and Lumley and Heagerty [17]. They developed a useful approach for the variance estimation of spatially correlated count data by combining the concept of window subsampling and the QL estimating function. (We discuss the contrast between their approach and ours below.)

We next review the concept of QL functions preliminary to our extension of QL functions to spatially dependent data. The definition of QL functions used in this paper is the same as that of McCullagh and Nelder ([20], page 327). On the assumption that the covariance matrix $\text{cov}(\mathbf{Y}) = \psi^2 \mathbf{V}$ is proportional to some function of the mean, the quasi-score function is given by

$$(1.1) \qquad \mathbf{U}(\boldsymbol{\beta}; \mathbf{Y}) = \mathbf{P}^T \mathbf{V}^{-1}[\mathbf{Y} - \boldsymbol{\theta}(\boldsymbol{\beta})]/\psi^2.$$

Here $\psi^2$ is a constant independent of $\boldsymbol{\theta}$ and $\mathbf{P}$ is the derivative matrix $\partial \boldsymbol{\theta}/\partial \boldsymbol{\beta}$. The quasi-score function (1.1) is still relevant when observations are dependent, although some restrictions must be imposed on the covariance matrix $\mathbf{V}$ in order that some scalar functions whose gradient vector is equal to $\mathbf{U}(\boldsymbol{\beta})$ exist ([20], page 333). At the observed $\mathbf{Y} = \mathbf{y}$, the QL estimate $\boldsymbol{\beta}$ is the solution of a QL estimating equation $\mathbf{U}(\hat{\boldsymbol{\beta}}; \mathbf{Y} = \mathbf{y}) = 0$. For nonlinear models, a closed-form expression for $\boldsymbol{\beta}$ is usually impossible to obtain. McCullagh [18] suggested investigating statistical properties of $\hat{\boldsymbol{\beta}}$ by considering $\hat{\boldsymbol{\beta}}$ as a root of the quasi-score function within a neighborhood of the true parameter point.

In the next section we follow this idea to establish consistency and asymptotic normality of QL estimates. Some assumptions are made such that the central limit theory for correlated response variables $\mathbf{Y}$ holds. An example is used to illustrate how QL estimating equations can test dependence between cross-classified response variables in Section 3.

**2. Asymptotic normality of quasi-likelihood estimates.** In this paper we focus on a logistic link function: $\mathrm{logit}(\boldsymbol{\theta}) = \mathbf{T}\boldsymbol{\beta}$, where $\boldsymbol{\beta} = (\beta_0, \ldots, \beta_u)^T$ and the $i$th row of $\mathbf{T}$ is $\mathbf{t}_i = (1, t_{i1}, \ldots, t_{iu})$ with $t_{ij} = 0$ or $1$. The individual parameters $\theta_i = E(Y(\mathbf{s}_i)|\mathbf{t}_i)$ are thus given by $\exp(\mathbf{t}_i^T\boldsymbol{\beta})/(1+\exp(\mathbf{t}_i^T\boldsymbol{\beta}))$ for $i = 1, \ldots, N$, and the $(ij)$th component of $\mathbf{P}$ of (1.1) is

$$(2.1) \qquad (\mathbf{P})_{ij} = \begin{cases} \theta_i(1-\theta_i), & \text{if } j = 1, \\ t_{i,j-1}\theta_i(1-\theta_i), & \text{if } j = 2, \ldots, u+1. \end{cases}$$

Let $\gamma_{ij} = \mathrm{corr}(Y(\mathbf{s}_i), Y(\mathbf{s}_j))$ denote the correlation between response variables at sites $\mathbf{s}_i$ and $\mathbf{s}_j$. In this paper $\mathbf{Y}$ is assumed to be an isotropic process, so that $\gamma_{ij}$ depends only on the distance between $\mathbf{s}_i$ and $\mathbf{s}_j$ and is independent of $\boldsymbol{\theta}$. For QL estimates to exist, the covariance structure must follow certain restrictions [26]. In our covariance matrix, the $(ij)$th entry of its inverse can be represented by $\mathbf{V}_{ij}^{-1} = (\sigma_i\sigma_j\gamma_{ij})^{-1}$, where $\sigma_i = \sqrt{\theta_i(1-\theta_i)}$. Thus $\mathbf{V}_{ij}^{-1}$ is independent of $\theta_k$ and $\partial \mathbf{V}_{ij}^{-1}/\partial \theta_k$, $\partial \mathbf{V}_{ik}^{-1}/\partial \theta_j$, $\partial \mathbf{V}_{kj}^{-1}/\partial \theta_i$ are equal to zero. This satisfies the condition for existence provided by McCullagh and Nelder ([20], page 334).

To ensure that the central limit theory holds, some restrictions should be imposed on $\gamma_{ij}$. One typical assumption is that the observations should satisfy long-range independence [6]. With this assumption, the random variables can be divided into blocks and treated as independent. The classical central limit theorem then leads to asymptotic normality for the sum of random variables immediately.

A condition more general than long-range independence is strong mixing [24]. For a stationary random field satisfying appropriate mixing conditions, Bolthausen [4] showed that the central limit theorem holds for grided data on $Z^k$. This version of a central limit theorem has been cited several times in the literature of spatial statistics [12, 13, 22]. Nevertheless, one of the assumptions for this central limit theorem is that the random fields are described with respect to an $L_\infty$ norm ([4], page 1047, and [22], page 56).

In contrast, in this paper we show the existence of asymptotic normality with respect to any $L_p$ norm by requiring correlations to decrease exponentially with distance. This allows us to extend the QL estimating function to the $L_1$ or Euclidean distance, which are the metrics most commonly used in spatial data. (Although we do not pursue the idea here, we conjecture that Assumption 2.1 could be adapted to more general conditions, including mixing conditions. Doing so would permit a tighter linkage between our results and those of the above-cited authors.)

ASSUMPTION 2.1. $\gamma_{ij} = a\rho^{d(\mathbf{s}_i,\mathbf{s}_j)}$, where $a$ is a positive constant such that $0 \leq a\rho \leq 1$, $\rho \in [0,1]$. $d(\mathbf{s}_{i_1}, \mathbf{s}_{i_2})$ denotes the $L_p$ distance between $\mathbf{s}_i$ and $\mathbf{s}_j$.

Let $c_k(\mathbf{s}_1,\ldots,\mathbf{s}_k)$ and $m_k(\mathbf{s}_1,\ldots,\mathbf{s}_k)$ denote the cumulant and product moment functions of order $k$, respectively, for centered variables $\mathbf{Z} = \mathbf{Y} - \boldsymbol{\theta}$. The Leonov–Shiryaev formula ([14], page 21) leads to

$$m_k(\mathbf{s}_{i_1},\ldots,\mathbf{s}_{i_k}) - \sum_{\Upsilon^{1,\ldots,k}_{k/2}} m_2(\mathbf{s}_{\omega_1},\mathbf{s}_{\omega_2}) \times \cdots \times m_2(\mathbf{s}_{\omega_{k-1}},\mathbf{s}_{\omega_k}) = c_k(\mathbf{s}_1,\ldots,\mathbf{s}_k)$$

for even $k$, where $\Upsilon^{1,\ldots,k}_{k/2}$ denotes the collection of all possible sets whose elements are $\frac{k}{2}$ disjoint pairs from $\{1,\ldots,k\}$. For odd $k$, $m_k(\mathbf{s}_1,\ldots,\mathbf{s}_k) = c_k(\mathbf{s}_1,\ldots,\mathbf{s}_k)$. We now make some assumptions on $c_k(\mathbf{s}_1,\ldots,\mathbf{s}_k)$.

ASSUMPTION 2.2. Some $\alpha_k < \frac{k}{2}$ exist such that $\sum_{i_1 \neq \cdots \neq i_k}^{N} c_k(\mathbf{s}_{i_1},\ldots,\mathbf{s}_{i_k}) = O(N^{\alpha_k})$.

There are a number of processes for which this assumption holds, including examples such as independent and $m$-dependent processes ([5], page 20). In a separate paper we hope to characterize more thoroughly the set of processes satisfying Assumption 2.2. We use Assumptions 2.1 and 2.2 to establish Theorem 2.1, whose proof is outlined in the Appendix.

THEOREM 2.1. *For binary variables* $\mathbf{Y}$ *satisfying Assumptions* 2.1 *and* 2.2, $\frac{1}{\sqrt{N}}\mathbf{a}'(\mathbf{Y} - \boldsymbol{\theta}) \sim N(\mathbf{0}, \frac{1}{N}\mathbf{a}'\mathbf{V}\mathbf{a}) + O_P(\frac{1}{\sqrt{N}})$ *for any bounded vector* $\mathbf{a}$.

From a geometrical perspective, the matrix $\mathbf{P}^T\mathbf{V}^{-1}$ in the QL estimating equations represents a projection matrix of the residual vector $\mathbf{y} - \boldsymbol{\theta}$ onto the space spanned by the columns of $\mathbf{T}$. Before showing consistency of the QL estimates, we first study some properties of this matrix.

NOTATION. For matrices $A$ and $B$, the symbol $A \succeq B$ means that, for any $i$ and $j$, the $(ij)$th entry of $A - B$ is nonnegative.

To simplify matters, we assume that the $m \times n$ lattice of locations $\mathbf{s}_i$ is labeled columnwise, so that the first column of locations is labeled $\mathbf{s}_1, \mathbf{s}_2, \ldots, \mathbf{s}_m$, the second column is labeled $\mathbf{s}_{m+1}, \mathbf{s}_{m+2}, \ldots, \mathbf{s}_{2m}$, and so on. Also, in the following proofs, the notation $\mathbf{\Gamma}_a(\rho; L_p)$ is used to indicate that the correlation matrix depends on $a$, $\rho$ and the $L_p$ distance. $a$ is assumed to lie within $(0, 1]$ because if $a = 0$, $\mathbf{\Gamma}$ is a zero matrix and there is nothing to discuss.

LEMMA 2.1. $\max\{|(\mathbf{\Gamma}_a^{-1}(\rho; L_1))_{ij}| : i, j = 1, \ldots, N\}$ involves only $a$ and $\rho$.

PROOF. First we study the correlation matrix $\mathbf{\Gamma}_a(\rho; L_1)$ for $a = 1$ and then extend the result to other $a \in (0, 1)$. Let $\mathbf{\Omega}_n$ denote an $n \times n$ matrix with $(\mathbf{\Omega}_n)_{ij} = 1$ if $i = j$ and $\rho^{|i-j|}$ otherwise. Then, $\mathbf{\Gamma}_1(\rho; L_1) = \mathbf{\Omega}_n \otimes \mathbf{\Omega}_m$, where $\otimes$ is the Kronecker product. It is possible to show ([23], page 255) that $\mathbf{\Gamma}_1^{-1}(\rho; L_1) = \mathbf{\Omega}_n^{-1} \otimes \mathbf{\Omega}_m^{-1}$, where

$$(2.2) \quad \mathbf{\Omega}_n^{-1} = (1 - \rho^2)^{-1} \begin{pmatrix} 1 & -\rho & 0 & 0 & \cdots & 0 & 0 \\ -\rho & 1+\rho^2 & -\rho & 0 & \cdots & 0 & 0 \\ \vdots & \vdots & \vdots & \vdots & & \vdots & \vdots \\ 0 & 0 & 0 & 0 & \cdots & -\rho & 1 \end{pmatrix}_{n \times n}.$$

Therefore, the maximum absolute value in $\mathbf{\Gamma}_1^{-1}(\rho; L_1)$, say $\rho^*$, only involves $a$ and $\rho$.

For a general matrix $\mathbf{\Gamma}_a(\rho)$, $a \in (0, 1)$, first notice that $\mathbf{\Gamma}_a(\rho; L_1)$ is a positive-definite matrix. Since $\mathbf{\Gamma}_1(\rho; L_1) \succeq \mathbf{\Gamma}_a(\rho; L_1) \succeq a\mathbf{\Gamma}_1(\rho; L_1)$, applying the strong partial ordering of positive-definite matrices ([20], page 335) gives

$$(2.3) \quad \frac{1}{a}\mathbf{\Gamma}_1^{-1}(\rho; L_1) \succeq \mathbf{\Gamma}_a^{-1}(\rho; L_1) \succeq \mathbf{\Gamma}_1^{-1}(\rho; L_1).$$

Thus it follows that the maximum value in $\mathbf{\Gamma}_a^{-1}(\rho; L_1)$ is bounded between $\rho^*/a$ and $\rho^*$, and these only involve $a$ and $\rho$. □

LEMMA 2.2. With correlations under the $L_1$ metric, $\max\{|(\mathbf{P}^T\mathbf{V}^{-1})_{ij}| : i, j = 1, \ldots, N\} \leq C_0$ for some constant $C_0$ independent of $N$. Moreover, some constant $C_1$ exists such that $\limsup_{N \to \infty} \frac{1}{N}(\mathbf{P}^T\mathbf{V}^{-1}\mathbf{P})_{kl} \leq C_1$ for all $k, l = 1, \ldots, N$.

PROOF. For convenience in the proof, let $\mathbf{V}_a$ denote the covariance matrix corresponding to $\mathbf{\Gamma}_a(\rho; L_1)$. To evaluate the first row of $\mathbf{P}^T\mathbf{V}_1^{-1}$, note from (2.2) that each row or column of $\mathbf{\Omega}_n^{-1}$ has at most three nonzero entries.

The Kronecker product then implies that each column of $\mathbf{\Gamma}_1^{-1}(\rho; L_1)$ has at most nine nonzero entries. So, simple algebra gives $|(\mathbf{P}'\mathbf{V}_1^{-1})_{1j}| \leq 9\rho^*\theta^*$, where $\rho^*$ is given in the proof of Lemma 2.1 and $\theta^* = \max\{\sqrt{\frac{\theta_l(1-\theta_l)}{\theta_k(1-\theta_k)}} : l, k = 1, \ldots, N\}$. This $\theta^*$ is independent of $N$ because at most $2^u$ possible values of $\theta$ exist due to $\theta_k$ only involving $\boldsymbol{\beta}$ and the binary vector $\mathbf{t}_k$. From (2.1), the distinction between the first and other rows of $\mathbf{P}^T\mathbf{V}_1^{-1}$ is the multiplier $t_{ij}$. It then follows that $|(\mathbf{P}^T\mathbf{V}_1^{-1})_{ij}| \leq |(\mathbf{P}^T\mathbf{V}_1^{-1})_{1j}|$ for all $i$ and $j$ because $t_{ij} = 0$ or $1$. This implies that $9\rho^*\theta^*$ is an upper bound for the absolute values in $\mathbf{P}^T\mathbf{V}_1^{-1}$.

Next, it is easy to see that $|(\mathbf{P}^T\mathbf{V}_1^{-1}\mathbf{P})_{kl}| \leq \sum_{m=1}^{N}|(\mathbf{P}^T\mathbf{V}_1^{-1})_{km}(\mathbf{P}^T)_{ml}| \leq \frac{9}{4}N\rho^*\theta^*$ for all $k$ and $l$ because all entries of $\mathbf{P}$ in (2.1) are not larger than $\frac{1}{4}$. Consequently, $\limsup_{N\to\infty} \frac{1}{N}(\mathbf{P}^T\mathbf{V}_1^{-1}\mathbf{P})_{kl} \leq C$ for some constant $C$.

Finally, to generalize the above to $V_a$, $a \in (0,1)$, note that all entries in $\mathbf{P}$ and $\boldsymbol{\Sigma}^{-1/2}$ are nonnegative. Thus (2.3) implies that

$$\text{(2.4)} \qquad \frac{1}{a}\mathbf{P}^T\mathbf{V}_1^{-1} \succeq \mathbf{P}^T\mathbf{V}^{-1} \succeq \mathbf{P}^T\mathbf{V}_1^{-1}$$

and

$$\text{(2.5)} \qquad \frac{1}{a}\mathbf{P}^T\mathbf{V}_1^{-1}\mathbf{P} \succeq \mathbf{P}^T\mathbf{V}^{-1}\mathbf{P} \succeq \mathbf{P}^T\mathbf{V}_1^{-1}\mathbf{P}.$$

Applying the results for $\mathbf{V}_1$ to (2.4) and (2.5) gives the desired result. □

The results of Lemma 2.2 can now be extended to any $L_p$ space, $1 \leq p < \infty$. The following theorem provides this general result.

THEOREM 2.2. *Lemma* 2.2 *holds for the $L_p$ metric, $1 \leq p < \infty$.*

PROOF. By Minkowski's inequality, we know that the $L_1$ metric is greater than the other $L_p$ metrics, $1 \leq p < \infty$. Since $\rho$ is between 0 and 1, we have $\mathbf{\Gamma}_a(\rho; L_p) \succeq \mathbf{\Gamma}_a(\rho; L_1)$ and thus the strong partial ordering of positive-definite matrices implies $(\mathbf{\Gamma}_a(\rho; L_1))^{-1} \succeq (\mathbf{\Gamma}_a(\rho; L_p))^{-1}$. Lemma 2.2 and an argument similar to the discussion of (2.4) and (2.5) give the desired result. □

The previous theorems show that $\mathbf{P}^T\mathbf{V}^{-1}$ is a bounded vector. So, by Theorem 2.1, $\mathbf{U}(\boldsymbol{\beta})$ is asymptotically normal. We now focus on the $L_1$ and $L_2$ distances. For simplicity, we note that the role of $\mathbf{P}^T\mathbf{V}^{-1}\mathbf{P}$ in quasi-score functions is similar to Fisher's information in ordinary likelihood functions, and therefore $\mathbf{I}(\boldsymbol{\beta})$ is used below to denote $\mathbf{P}^T\mathbf{V}^{-1}\mathbf{P}$.

THEOREM 2.3. *For $\mathbf{Y}$ satisfying Assumptions* 2.1 *and* 2.2, *we have*

$$\frac{1}{\sqrt{N}}\mathbf{U}(\boldsymbol{\beta}) \sim N\left(\mathbf{0}, \frac{1}{N}\mathbf{I}(\boldsymbol{\beta})\right) + O_P\left(\frac{1}{\sqrt{N}}\right) \qquad \text{as } N \to \infty.$$

Next we derive the limiting distribution of the QL estimate.

NOTATION. (a) We use $\mathbf{O}(N^p)$ to denote a matrix $A$ satisfying $\limsup_{N\to\infty} \frac{1}{N^p}(A)_{ij} \leq K$ for all $i$ and $j$ and for some constant $K$. (b) abs($A$) represents a matrix whose $(ij)$th element equals the absolute value of the $(ij)$th element of $A$. (c) $\mathbf{D}_{\beta_j}$ and $\mathbf{D}_\beta$ denote partial derivatives with respect to $\beta_j$ and $\boldsymbol{\beta}$, respectively.

ASSUMPTION 2.3. *The matrix $\frac{1}{N}\mathbf{D}_\beta\mathbf{U}(\boldsymbol{\beta})$ is negative definite at the true parameter $\boldsymbol{\beta}_0$ with probability going to 1 as $N\to\infty$.*

Assumption 2.3 is a common assumption used for the derivative of score functions (e.g., [19]).

LEMMA 2.3. *The QL estimate $\hat{\boldsymbol{\beta}}$ is consistent.*

PROOF. It is easy to see that

$$\mathbf{D}_\beta\mathbf{U}(\boldsymbol{\beta}) = (\mathbf{D}_\beta\mathbf{P}^T)\mathbf{V}^{-1}(\mathbf{Y}-\boldsymbol{\theta}) + \mathbf{P}^T(\mathbf{D}_\beta\mathbf{V}^{-1})(\mathbf{Y}-\boldsymbol{\theta}) - \mathbf{I}(\boldsymbol{\beta}).$$

Let $(t_{10},\ldots,t_{N0})^T = \mathbf{1}$ represent the first column of the design matrix and let $\boldsymbol{\xi}_i = (t_{1i}\theta_1(1-\theta_1),\ldots,t_{Ni}\theta_N(1-\theta_N))$, $i=0,\ldots,u$, represent the $(i+1)$st row of $\mathbf{P}^T$. Then $\mathbf{D}_{\beta_j}\boldsymbol{\xi}_i = (t_{1i}t_{1j}\theta_1(1-\theta_1)(1-2\theta_1),\ldots,t_{Ni}t_{Nj}\theta_N(1-\theta_N)(1-2\theta_N))$ for $j=0,\ldots,u$. Since the $t_{ij}$ are binary and $|1-2\theta_i|\leq 1$, it follows that $\boldsymbol{\xi}_i \succeq \mathrm{abs}(\mathbf{D}_{\beta_j}\boldsymbol{\xi}_i)$ and thus $\mathbf{P}^T \succeq \mathrm{abs}(\mathbf{D}_\beta\mathbf{P}^T)$. Similarly, we can show that $\boldsymbol{\Sigma}^{-1/2} \succeq \mathrm{abs}(\mathbf{D}_{\beta_j}\boldsymbol{\Sigma}^{-1/2})$, where $\boldsymbol{\Sigma} = \mathrm{diag}(\theta_i(1-\theta_i))$, and thus $2\mathbf{V}^{-1} \succeq \mathrm{abs}(\mathbf{D}_\beta\mathbf{V}^{-1})$.

Therefore, an argument similar to the proof of Theorem 2.2 gives that $(\mathbf{D}_\beta\mathbf{P}^T)\mathbf{V}^{-1}$ and $\mathbf{P}^T(\mathbf{D}_\beta\mathbf{V}^{-1})$ are $\mathbf{O}(1)$. By Theorem 2.1, a normal variable $Z^*$ exists for any bounded vector $\mathbf{a}$ such that $\frac{1}{\sqrt{N}}\mathbf{a}^T(\mathbf{Y}-\boldsymbol{\theta}) = Z^* + \mathbf{O}_P(\frac{1}{\sqrt{N}})$ as $N\to\infty$. It then follows from the Cramér–Wold device and Theorem 2.1 that

$$(\mathbf{D}_\beta\mathbf{P}^T)\mathbf{V}^{-1}(\mathbf{Y}-\boldsymbol{\theta}) + \mathbf{P}^T(\mathbf{D}_\beta\mathbf{V}^{-1})(\mathbf{Y}-\boldsymbol{\theta}) = \mathbf{O}_P(\sqrt{N})$$

and therefore

(2.6) $\quad\quad \mathbf{D}_\beta\mathbf{U}(\boldsymbol{\beta}) = \mathbf{O}_P(\sqrt{N}) - \mathbf{I}(\boldsymbol{\beta}) \quad$ as $N\to\infty$.

As a result, $\frac{1}{N}\mathbf{D}_\beta\mathbf{U}(\boldsymbol{\beta})|_{\boldsymbol{\beta}=\boldsymbol{\beta}_0} \to -\frac{1}{N}\mathbf{I}(\boldsymbol{\beta}_0)$ with probability going to 1 as $N\to\infty$. From Assumption 2.4, therefore $\frac{1}{N}\mathbf{I}(\boldsymbol{\beta}_0)$ has a positive-definite limit.

It follows next by the inverse theorem [3] that an open ball $B(\boldsymbol{\beta}_0,r)$ exists such that $\frac{1}{N}\mathbf{U}(\boldsymbol{\beta})$ is one-to-one on the ball with probability going to 1. Also, $\frac{1}{N}\mathbf{U}(B(\boldsymbol{\beta}_0,r))$ contains an open ball $B(\frac{1}{N}\mathbf{U}(\boldsymbol{\beta}_0),r^*)$ for some $r^*$.

By Theorem 2.3, $\frac{1}{N}\mathbf{U}(\boldsymbol{\beta}_0) \to E(\mathbf{U}(\boldsymbol{\beta}_0)) = \mathbf{0}$ in probability. Hence for this $r^*$, $\|\frac{1}{N}\mathbf{U}(\boldsymbol{\beta}_0) - \mathbf{0}\| < r^*$ with probability going to 1, where $\|\cdot\|$ denotes the Euclidean norm. This then implies that $\mathbf{0} \in B(\frac{1}{N}\mathbf{U}(\boldsymbol{\beta}_0), r^*) \subseteq \frac{1}{N}\mathbf{U}(B(\boldsymbol{\beta}_0, r))$ with probability going to 1. Since $\frac{1}{N}\mathbf{U}(\boldsymbol{\beta})$ is one-to-one on $B(\boldsymbol{\beta}_0, r)$ and $r$ can be arbitrarily small, $\frac{1}{N}\mathbf{U}(\hat{\boldsymbol{\beta}}) = \mathbf{0}$ a.e. implies that $\hat{\boldsymbol{\beta}} \to \boldsymbol{\beta}_0$ in probability. □

THEOREM 2.4. *For $\mathbf{Y}$ satisfying Assumptions 2.1–2.3, the QL estimate $\hat{\boldsymbol{\beta}}$ has the limiting distribution*

$$\sqrt{N}(\hat{\boldsymbol{\beta}} - \boldsymbol{\beta}_0) \sim N(\mathbf{0}, N\mathbf{I}^{-1}(\boldsymbol{\beta}_0)) + \mathbf{O}_P\left(\frac{1}{\sqrt{N}}\right) \qquad as\ N \to \infty.$$

PROOF. The first-order Taylor series expansion gives

$$(2.7) \quad \mathbf{U}(\hat{\boldsymbol{\beta}}) = \mathbf{U}(\boldsymbol{\beta}_0) + \mathbf{D}_{\boldsymbol{\beta}}\mathbf{U}(\boldsymbol{\beta})|_{\boldsymbol{\beta}=\boldsymbol{\beta}_0}(\hat{\boldsymbol{\beta}} - \boldsymbol{\beta}_0) + \mathbf{o}_p(\|\hat{\boldsymbol{\beta}} - \boldsymbol{\beta}_0\|).$$

Assume that the inverse of $\mathbf{D}_{\boldsymbol{\beta}}\mathbf{U}(\boldsymbol{\beta})$ exists at $\boldsymbol{\beta}_0$. Then (2.7) implies that $(\hat{\boldsymbol{\beta}} - \boldsymbol{\beta}_0)(1 + \mathbf{o}_P(1)) = -(\mathbf{D}_{\boldsymbol{\beta}}\mathbf{U}(\boldsymbol{\beta}))^{-1}|_{\boldsymbol{\beta}=\boldsymbol{\beta}_0}\mathbf{U}(\boldsymbol{\beta}_0)$ because $\mathbf{U}(\hat{\boldsymbol{\beta}}) = \mathbf{0}$. In addition, it follows from (2.6) that

$$(\mathbf{D}_{\boldsymbol{\beta}}\mathbf{U}(\boldsymbol{\beta}))^{-1} = -(\mathbf{O}_P(\sqrt{N}) - \mathbf{I}(\boldsymbol{\beta}))^{-1}$$
$$= -(\mathbf{I}(\boldsymbol{\beta}))^{-1} + (\mathbf{I}(\boldsymbol{\beta}))^{-2} \cdot \mathbf{O}_P(\sqrt{N}) + \mathbf{O}\left(\frac{1}{N^2}\right).$$

Thus $(\hat{\boldsymbol{\beta}} - \boldsymbol{\beta}_0)(1 + o(1))$ can be written as

$$(\mathbf{I}(\boldsymbol{\beta}_0))^{-1}\mathbf{U}(\boldsymbol{\beta}_0) + (\mathbf{I}(\boldsymbol{\beta}_0))^{-2}\mathbf{O}_P\left(\frac{1}{\sqrt{N}}\right)\mathbf{U}(\boldsymbol{\beta}_0) + \mathbf{O}\left(\frac{1}{N^2}\right)\mathbf{U}(\boldsymbol{\beta}_0).$$

Because $\mathbf{I}(\boldsymbol{\beta}_0) = \mathbf{O}(N)$ and $\mathbf{U}(\boldsymbol{\beta}_0) = \mathbf{O}_p(\sqrt{N})$ from Theorems 2.1 and 2.2, it follows that $(\hat{\boldsymbol{\beta}} - \boldsymbol{\beta}_0)(1 + o_P(1)) = (\mathbf{I}(\boldsymbol{\beta}_0))^{-1}\mathbf{U}(\boldsymbol{\beta}_0) + \mathbf{O}_P(\frac{1}{N})$. Since $\hat{\boldsymbol{\beta}}$ is consistent, this implies that $\hat{\boldsymbol{\beta}} - \boldsymbol{\beta}_0 = (\mathbf{I}(\boldsymbol{\beta}_0))^{-1}\mathbf{U}(\boldsymbol{\beta}_0) + \mathbf{O}_P(\frac{1}{N})$ as $N \to \infty$. The theorem then follows immediately from Theorem 2.3. □

The average rate of convergence for $\hat{\boldsymbol{\beta}}$ to $\boldsymbol{\beta}_0$ can be computed as follows.

COROLLARY 2.1. $E(\hat{\boldsymbol{\beta}} - \boldsymbol{\beta}_0) = \mathbf{O}(\frac{1}{N})$ *as $N \to \infty$.*

PROOF. This result follows immediately from Theorem 2.4. □

It follows from Theorem 2.4 that $E(\hat{\boldsymbol{\beta}}) \to \boldsymbol{\beta}_0$ and $\mathrm{cov}(\boldsymbol{\beta}) \to \mathbf{I}^{-1}(\boldsymbol{\beta}_0)$ as the sample size increases. In the terminology of Godambe [10], it can be shown that the QL estimating equation is an asymptotically unbiased optimal estimating equation [19].

**3. Application of quasi-likelihood functions.** In this section we illustrate the use of the previously developed theory to data cited by Fingleton [9]. These well-known data [8, 9] were collected from Lansing Woods to examine whether the presence of hickory near a sample site tends to discourage the presence of maple.

For notational convenience, let $X(\mathbf{s})$ and $Y(\mathbf{s})$ denote the corresponding indicator variables for the presence of hickory and maple at site $\mathbf{s}$, respectively. We define the leftmost of a horizontal pair as the "first" and the uppermost of a vertical pair as the "first" as Fingleton did in his 1986 paper. One approach for testing independence between two correlated binary variables is to model these two variables by a conditional logistic model,

$$(3.1) \quad \theta_i = P(Y(\mathbf{s}_i) = 1 | X(\mathbf{s}_i) = x_i) = \frac{\exp(\beta_0 + \beta_1 x_i)}{1 + \exp(\beta_0 + \beta_1 x_i)}, \qquad i = 1, \ldots, 256.$$

That is, at a specific site $\mathbf{s}_i$, given the information whether hickory is present or not, the model indicates how likely it is that maple is present. Thus, if the estimate of $\beta_1$ in (3.1) is not significantly different from zero, we may say that $\mathbf{X}$ and $\mathbf{Y}$ are independent. The concept of this model is not very far from the traditional chi-squared approach. In fact, when no spatial dependence exists among observations, this approach is asymptotically equivalent to a chi-squared test in the analysis of contingency tables [1].

In geostatistics the elements of the matrix of correlations between sites are usually obtained from semivariogram models parameterized by constants denoting the nugget effect, the sill and the range [7]. For the maple data, we fit the exponential semivariogram model depicted in Figure 1. The correlation of maples at sites $\mathbf{s}_i$ and $\mathbf{s}_j$ can therefore be estimated by $\rho_y(\mathbf{s}_i, \mathbf{s}_j) = \exp(-d(\mathbf{s}_i, \mathbf{s}_j)/1.091)$, where 1.091 comes from the "effective" range of the fitted exponential variogram model and $d(\mathbf{s}_i, \mathbf{s}_j)$ is the $L_2$ distance between sites $\mathbf{s}_i$ and $\mathbf{s}_j$.

Since the quasi-score function is nonlinear, the Newton–Raphson method, $\hat{\boldsymbol{\beta}}_{j+1} = \hat{\boldsymbol{\beta}}_j + (\hat{\mathbf{P}}'_j \hat{\mathbf{V}}_j^{-1} \hat{\mathbf{P}}'_j)^{-1} \hat{\mathbf{P}}'_j \hat{\mathbf{V}}_j^{-1} (\mathbf{y} - \hat{\boldsymbol{\theta}}_j)$, was used to derive a numerical solution. An approximation of $\hat{\boldsymbol{\beta}}$ obtained by iteration is $(\hat{\beta}_0, \hat{\beta}_1) = (-0.34, -0.26)$ and $\text{vâr}(\hat{\beta}_0) = 0.138$, $\text{vâr}(\hat{\beta}_1) = 0.001$ and $\text{côv}(\hat{\beta}_0, \hat{\beta}_1) = -0.003$.

According to Theorem 2.4, we know that $\hat{\boldsymbol{\beta}}$ is asymptotically normal. So, we can construct an approximate chi-squared test for the hypothesis that the parameter of independence $\beta_1$ is zero by using $\hat{\beta}_1^2 / \hat{\sigma}^2(\hat{\beta}_1) = 6.76$. This value is between those of the traditional chi-squared test (25.7) and Fingleton's deflated chi-squared value of 2.77 [9]. Therefore we may say that our approach provides a balanced point of view between the conservative traditional chi-squared test and the liberal Fingleton deflated chi-squared test. In addition, the negative $\hat{\beta}_1$ ($-0.26$) with a significant $p$-value (0.007 comparing to the $\chi_1^2$ table) shows strong evidence that the presence of hickory discourages the presence of maple. Note that neither the traditional chi-squared test nor Fingleton's method could tell us whether the interaction between maple and hickory is positive or negative.

**4. Discussion.** Analysis of non-Gaussian spatially correlated data is usually difficult because of the complexity of the associated distributional forms. In this paper we have extended quasi-likelihood estimating equations to deal with spatially correlated data when the response variables are binary. We model the marginal response probability by a logistic regression; one benefit of this is that we can avoid the specification of likelihood functions. To employ the proposed method we only need to know the means and covariances between observations. If correlations between observations are decreasing exponentially with distance (a reasonable assumption in many spatial settings), the QL estimates are shown to be asymptotically normal and consistent. This allows us to estimate and test fixed effects for data with correlation between sites.

Although the QL approach discussed in this paper was originally motivated by a question in spatial statistics, the proposed approach is potentially applicable in other settings. In previous literature, estimating equation approaches developed for dependent observations have mostly focused on repeated measurements with block-diagonal covariance matrices, that is, with independence between subjects. The QL approach developed in this paper generalizes such methods; in this paper we show that, under good control on correlations, the asymptotics of QL estimates can hold for non-block-diagonal covariance matrices. This opens the possibility of applying QL ap-proaches to other situations where the data do not exhibit independence between subjects.

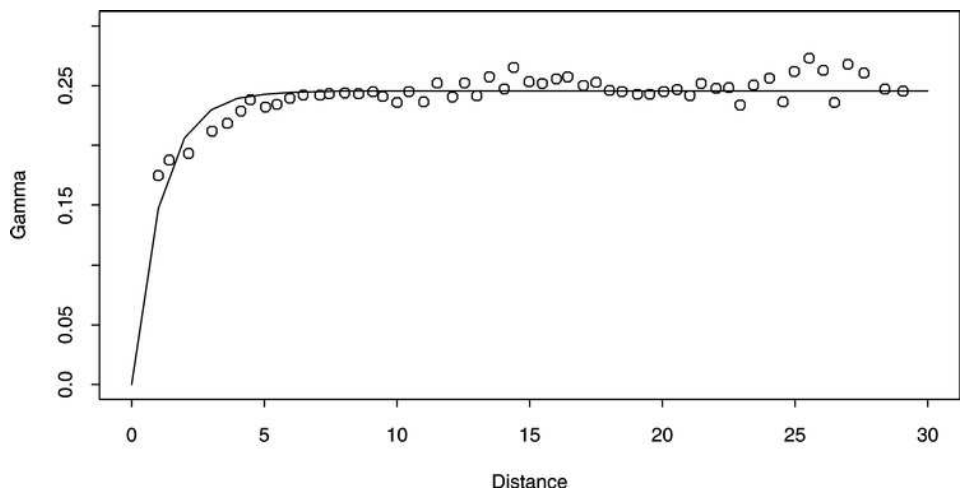

FIG. 1. *Empirical semivariogram for the maple data. The solid line indicates the fitted variogram model with absolute sill* 0.246 *and range* 1.091.

We briefly outline a number of extensions of this work which we are currently pursuing. First, the proposed approach developed in this paper is for complete data. The forest ecology problem referred to in Section 1 is not presented in this paper partly because observations are missing from several of the sites. Interestingly, the problem of missing data is related to the problem of analyzing data on an irregular lattice, another issue of note. In addition, the forest ecology data show possible large-scale spatial trends among observations. Conceivably this could be addressed by including appropriate covariates in the modeling work, although it remains to be seen how this would influence the asymptotics. On the other hand, as suggested by an anonymous referee, another possible solution for this problem is to fit two estimating equations simultaneously, one for the mean structure and the other for correlation structure. This approach was taken by Zeger [26] and McShane, Albert and Palmatier [21] for time series data by conditioning on a latent process.

Our work in this paper has been restricted to exponential variograms, although the flexibility exists to allow different $L_p$ metrics. Nonetheless, it would be interesting to extend the methods of this paper to other variogram models that frequently occur in the analysis of spatial data. Moreover, the proposed method in this paper could be generalized to spatially correlated count data. However, this would involve reconstructing the covariance structure, which itself would be a considerable task. We leave this to future work.

## APPENDIX

Here we list some important results used in the proof of Theorem 2.1. $\mathbf{Z}$ is used to denote $\mathbf{Y} - \boldsymbol{\theta}$.

LEMMA A.1. *For centered binary variables $Z_i = Y_i - \theta_i$, we have*

$$E[Z_1^2 Z_2^{p_2} \cdots Z_k^{p_k}] = \mathrm{Var}(Z_1) E[Z_2^{p_2} \cdots Z_k^{p_k}] + (1 - 2\theta_1) E[Z_1 Z_2^{p_2} \cdots Z_k^{p_k}].$$

PROOF. The expectation $E[Z_1^2 Z_2^{p_2} \cdots Z_k^{p_k}]$ can be expanded to

$$\begin{aligned}
E[Z_1^2 &Z_2^{p_2} \cdots Z_k^{p_k}] \\
&= \theta_1^2 E[Z_2^{p_2} \cdots Z_k^{p_k}] \\
&\quad + (1 - 2\theta_1) \sum_{(z_2,\ldots,z_k)} \cdots \sum z_2^{p_2} \cdots z_k^{p_k} \\
&\qquad\qquad\qquad\qquad \times P(Z_1 = 1 - \theta_1, Z_2 = z_2, \ldots, Z_k = z_k).
\end{aligned} \quad \text{(A.1)}$$

The $k$-fold summation in the second line can be converted to $E(Z_1 Z_2^{p_2} \cdots Z_k^{p_k}) + \theta_1 E(Z_2^{p_2} \cdots Z_k^{p_k})$. Inserting this result back into (A.1) gives the desired result. □

Lemma A.1 focuses on the "dependent" item in the expectation. When $Z_1$ is independent of $(Z_2, \ldots, Z_k)$, $E(Z_1^2 Z_2^{p_2} \cdots Z_k^{p_k}) = \text{var}(Z_1) E(Z_2^{p_2} \cdots Z_k^{p_k})$, which is exactly the same form of Lemma A.1 when $\theta_1 = 0.5$. Thus $(1 - 2\theta_1) E(Z_1 Z_2^{p_2} \cdots Z_k^{p_k})$ of Lemma A.1 can be considered to be the impact of dependence on the expected value, and we can expect that this impact is reduced for $\theta_1$ close to 0.5. In fact, $E(Z_1 Z_2 Z_3) = 0$ if $Z_1, Z_2$ and $Z_3$ are from a truncated Gaussian random field.

LEMMA A.2. *For variables $\mathbf{Z}$ satisfying Assumption 2.1 with the $L_1$ metric,*

$$\sum_{i_1 \neq i_2}^{N} \sum^{N} \text{cov}(Z(\mathbf{s}_{i_1}), Z(\mathbf{s}_{i_2})) \leq \frac{2\rho - \rho^2}{(1-\rho)^2} aN.$$

PROOF. The sum of correlations between a given site $\mathbf{s}_{i_1}$ and all other sites over a region has a maximum value when $\mathbf{s}_{i_1}$ is the center. So $\sum_{i_2=1}^{N} \text{cov}(X(\mathbf{s}_{i_1}), X(\mathbf{s}_{i_2})) \leq \sum_{\mathbf{s}_{i_2} \in \mathcal{R}_m} \rho^{d(([bt]\lfloor m/2 \rfloor+1, [bt]\lfloor n/2 \rfloor+1), \mathbf{s}_{i_2})}$, where $\mathcal{R}_m$ denotes the upper right quarter of the region. This holds because $\text{var}(Z(\mathbf{s})) \leq \frac{1}{4}$ and the sum of correlations in each quarter is the same. This lemma then follows from

$$\sum_{\mathbf{s}_{i_2} \in \mathcal{R}_m} \rho^{d(([bt]\lfloor m/2 \rfloor+1, [bt]\lfloor m/2 \rfloor+1), \mathbf{s}_{i_2})} = -1 + \sum_{u=0}^{m/2} \sum_{t=0}^{n/2} \rho^{t+u} \leq \frac{2\rho - \rho^2}{(1-\rho)^2}. \quad \square$$

LEMMA A.3. *For variables $\mathbf{Z}$ satisfying Assumption 2.1 with the $L_2$ metric,*

$$\sum_{i_1 \neq i_2}^{N} \sum^{N} \text{cov}(Z(\mathbf{s}_{i_1}), Z(\mathbf{s}_{i_2})) \leq \left( \frac{2\rho}{1-\rho} + \frac{\pi}{2} \left( \frac{1}{\log(\rho)} \right)^2 \right) aN.$$

PROOF. The sum of correlations between the center of the study region and all sites in $\mathcal{R}_m$ is $-1 + \sum_{u=0}^{m/2} \sum_{t=0}^{n/2} \rho^{\sqrt{t^2+u^2}}$, which can be shown to be less than $\frac{2\rho}{1-\rho} + \frac{\pi}{2}(\frac{1}{\log(\rho)})^2$ by an integral test. An argument similar to Lemma A.2 gives the desired result. □

Lemmas A.2 and A.3 are used to control correlations. Other details for controlling higher-order correlations are shown in [16] and omitted here.

Let $\phi(s)$ and $\varphi(s)$ denote the corresponding characteristic functions of $\frac{1}{\sqrt{N}} \mathbf{a}' \mathbf{Z}$ and a normal random variable with mean zero and variance $\frac{1}{N} \mathbf{a}' \mathbf{V} \mathbf{a}$, respectively. Our approach to show Theorem 2.1 is to prove that $\phi(s) = \varphi(s) + O(\frac{1}{\sqrt{N}})$. Let $\phi^*(s)$ denote the $N$th truncated Taylor series of $\phi(s)$. Under Assumptions 2.1 and 2.2, we can show that the imaginary part of $\phi^*(s)$ is $O(\frac{1}{\sqrt{N}})$ and the real part $\phi^*(s)$ has an approximation

$$
(A.2) \quad 1 + \sum_{q=1}^{\lfloor N/2 \rfloor} \frac{s^{2q}}{q! 2^q} \frac{1}{N^q} \times \sum_{j=0}^{q} \binom{q}{j} \sum_{i_1 \neq \cdots \neq i_{q+j}}^{N} \cdots \sum^{N} \nu_{i_1} \cdots \nu_{i_{q-j}} \xi_{i_{q-j+1}, i_{q-j+2}} \cdots \xi_{i_{q+j-1}, i_{q+j}}
$$

with error $O(\frac{1}{\sqrt{N}})$, where

$$\nu_i = a_i^2 \theta_i (1 - \theta_i) \quad \text{and} \quad \xi_{i,j} = a_i a_j [bt] \sqrt{\theta_i (1 - \theta_i) \theta_j (1 - \theta_j)} a \rho^{d(\mathbf{s}_i, \mathbf{s}_j)}.$$

The proof requires numerous steps involving combinatorial analysis. Readers interested in the process can find the details in [16].

In addition, it can be shown that

$$|\phi(s) - \phi^*(s)| \leq \frac{|s|^N}{(N+1)!} E \min\left( |s| \cdot \left| \frac{\mathbf{a}' \mathbf{Z}}{\sqrt{N}} \right|^{N+1}, 2(N+1) \left| \frac{\mathbf{a}' \mathbf{Z}}{\sqrt{N}} \right|^N \right) = o(s^N).$$

The $N$th truncated Taylor series of $\varphi(s)$, denoted by $\varphi^*(s)$, can also be shown to have an approximation (A.2) with error $O(\frac{1}{N})$. The finite value of $\frac{1}{N} \mathbf{a}' \mathbf{V} \mathbf{a}$ from Lemma A.1 or A.2 (depending on the metric employed) then implies $|\varphi(s) - \varphi^*(s)| = o(s^N)$. Theorem 2.1 is then obvious from these results.

**Acknowledgment.** The authors are grateful to an anonymous referee for helpful suggestions and for making us aware of a broader range of related literature.

DEPARTMENT OF MATHEMATICS
NATIONAL CHUNG CHENG UNIVERSITY
160 SAN-HSING, MIN-HSIUNG
CHIA-YI 621
TAIWAN
REPUBLIC OF CHINA
E-MAIL: pslin@math.ccu.edu.tw

DEPARTMENT OF STATISTICS
UNIVERSITY OF WISCONSIN–MADISON
1300 UNIVERSITY AVENUE
MADISON, WISCONSIN 53706
USA
E-MAIL: clayton@stat.wisc.edu